\documentclass[11pt]{article} \usepackage{amssymb}
\usepackage{amsfonts} \usepackage{amsmath} \usepackage{amsthm} \usepackage{bm}
\usepackage{latexsym} \usepackage{epsfig}

\newtheorem*{theorem*}{Theorem}
\newtheorem{theorem}{Theorem}[section]
\newtheorem{lemma}[theorem]{Lemma}
\newtheorem*{proposition*}{Proposition}

\newtheorem{proposition}[theorem]{Proposition}

\newtheorem{conjecture}[theorem]{Conjecture}

\newcommand{\ignore}[1]{}

\newcommand{\enote}[1]{} \newcommand{\knote}[1]{}
\newcommand{\rnote}[1]{}

\newcommand{\E}{{\bf E}} \newcommand{\Cov}{{\bf Cov}}
\newcommand{\Var}{{\bf Var}}

\newcommand{\N}{\mathbb N} \newcommand{\R}{\mathbb R}

\newcommand{\half}{{\textstyle \frac12}}

\newcommand{\fourth}{{\textstyle \frac14}}
\newcommand{\fifth}{{\textstyle \frac15}}

\renewcommand{\phi}{\varphi}

\begin{document}
\title{Iterative Maximum Likelihood on Networks}
\author{
  Elchanan Mossel\footnote{Weizmann Institute and U.C. Berkeley. E-mail:
    mossel@stat.berkeley.edu. Supported by a Sloan fellowship in
    Mathematics, by BSF grant 2004105, by NSF Career Award (DMS 054829)
    by ONR award N00014-07-1-0506 and by ISF grant 1300/08}~ and
  Omer Tamuz\footnote{9 Smilanski St., Herzlia 46361, Israel}
 }

\date{\today}
\maketitle
\begin{abstract}
  We consider $n$ agents located on the vertices of a connected graph.
  Each agent $v$ receives a signal $X_v(0)\sim N(\mu,1)$
  where $\mu$ is an unknown quantity. A natural iterative
  way of estimating $\mu$ is to perform the following
  procedure. At iteration $t+1$ let $X_v(t+1)$
  be the average of $X_v(t)$ and of $X_w(t)$ among all the
  neighbors $w$ of $v$. It is well known that this procedure
  converges to $X(\infty) = \frac{1}{2} |E|^{-1} \sum d_v X_v$ where
  $d_v$ is the degree of $v$.

  In this paper we consider a variant of simple iterative averaging,
  which models ``greedy'' behavior of the agents.
  At iteration $t$, each agent $v$ declares the value of its estimator
  $X_v(t)$ to all of its neighbors. Then, it updates $X_v(t+1)$
  by taking the maximum likelihood (or minimum variance) estimator of
  $\mu$, given $X_v(t)$ and
  $X_w(t)$ for all neighbors $w$ of $v$, and the structure of the graph.

  We give an explicit efficient procedure for calculating $X_v(t)$,
  study the convergence of the process as $t \to \infty$ and show that if
  the limit exists then $X_v(\infty) = X_w(\infty)$ for all $v$ and $w$.
 For graphs that are
  symmetric under actions of transitive groups, we show that
  the process is efficient.  Finally, we show that the greedy
  process is in some cases
  more efficient than simple averaging, while in other cases
  the converse is true, so that, in this model, ``greed'' of the
  individual agents may or may not have an adverse
  affect on the outcome.

  The model discussed here may be viewed as the Maximum-Likelihood
  version of models studied in Bayesian Economics.  The ML variant is
  more accessible and allows in particular to show the significance of
  symmetry in the efficiency of estimators using networks of agents.

\end{abstract}

\newpage
\section{Introduction}
Networks and graphs are often viewed as computational models. In
computational complexity several complexity classes are studied in
terms of corresponding computation graphs, for example
finite-automata, PSPACE and
LOG-SPACE. For general background see, e.g.,~\cite{ComplexityZOO}.
In parallel computing, networks are used to model the
communication network between different computers,
while in sparse sensing the connectivity network is of fundamental computation
significance (see, e.g., \cite{AmmariDas:08} and \cite{ZhangHou:08}).

A recent trend emanating from Economics and Game Theory considers
networks where different nodes correspond to computational entities
with different objectives \cite{Jackson:06}.
Recent models in Bayesian Economics
consider models where each player is repeatedly taking actions that
are based on a signal he has received that is correlated with the
state of the word and past actions of his neighbors
(\cite{BikhHirshWelch:08}, \cite{Banerjee:08}, \cite{SmithSorensen:00}).

In this paper we study a simple model where, in each iteration,
 agents iteratively try to optimally estimate the state of the
world, which is a single parameter $\mu \in \R$.
It is assumed that originally each agent receives an independent
sample from a normal distribution with mean $\mu$. Later at each
iteration each agent updates his estimate by taking the maximum
likelihood estimator of $\mu$ given its current estimator and those
of its neighbors, {\em and
  given the graph structure}. Note that for normal distributions,
the maximum likelihood estimator is identical with the minimum
variance unbiased estimator. At the first iteration, the
estimator at each node will be the average of the original signal at
the node and its neighbors.  However, from the second iteration on,
the procedure will not proceed by simple averaging due to the
correlation between the estimators at adjacent nodes. As we show
below, this correlation can be calculated given the structure of the
graph and results in dramatic differences from the simple averaging
process. Note that under this model, the agents are memoryless and
use only the results of the last iteration to calculate those the next.

The model suggested above raises a few basic questions:
\begin{itemize}
\item Is the process above well defined?
\item Can the estimators be efficiently calculated? Note that in the
Bayesian economic models (such as \cite{SmithSorensen:00}) there are
no efficient algorithms for updating beliefs.
\item Does the process converge?
\end{itemize}
We answer the first two questions positively, and conjecture that the answer
to the third is positive as well. Once these questions are addressed
we prove a number of results regarding the limit estimators including:
\begin{itemize}
\item We show that for connected graphs,
  as $t\to \infty$, the correlation between the estimators of the different agents
  goes to one.
\item We describe a graph for which the maximum likelihood
  process converges to an estimator different than the optimal.
\item We compare the statistical efficiency of the limiting estimator
  to the limiting estimator obtained by simple iterative averaging and
  to the optimal estimator, in different graphs.
\end{itemize}

\subsection{Formal Definition of the Model}

We consider a finite, undirected, connected graph
$G=(V,E)$, where each vertex has a self-loop so that $\forall v:\: (v,v)\in E$,
and a state of the world $\mu \in\R$.
We assign each vertex $v$ a normal unbiased estimator $X_v=X_v(0)$ of $\mu$ so that
$\E[X_v] = \mu$ and $\Var[X_v]=1$, for all $v$.
These estimators are uncorrelated.

In iteration $t\in\N$ we define
$X_v(t+1)$ to be the minimum variance unbiased estimator constructible over
the estimators of $v$ and its neighbors $N(v) = \{w | (v,w)\in E\}$ at time $t$
\begin{equation}
  \label{eq:model_def1}
X_v(t+1) = \sum_{w \in N(v)} \alpha_w X_w(t),\quad \mbox{ where: }
\end{equation}
\begin{equation}
  \label{eq:model_def2}
\sum_{w \in N(v)} \alpha_w = 1, \mbox{ and } \alpha \mbox{ minimizes } \Var[\sum_w \alpha_w X_w(t)].
\end{equation}
(note that $\alpha$ may be positive or negative).
The process $Y_v$ is given by simple iterative averaging so $Y_v(0):= X_v$, and
\begin{equation}
  \label{eq:y_def}
  Y_v(t+1)={1\over \mathrm{d_v}}\sum_{w\in N(v)}{Y_w(t)}.
\end{equation}
It is well known that $Y_v(t)$ converges to $Y_v(\infty)= \frac{1}{2} |E|^{-1} \sum d_v X_v$.

Finally, we define $Z(\infty)$ to be the global minimum variance unbiased
 estimator:
 \[
Z(\infty) = \frac{1}{|V|} \sum_{v \in V} X_v.
\]
Note that in the case of normally distributed $X_v$'s, the minimum
variance definitions coincide with those of maximum likelihood.

This scheme can be generalized to the case where the
original estimators $X_v$ have a general covariance structure,
with the definitions for
$X_v(t)$ and $Y_v(t)$ remaining essentially the same, and that
of $Z(\infty)$ changing to the form of Eq.~\ref{eq:algorithm} below.

\subsection{Statements of the main results}
\begin{itemize}
\item The process defined by Eqs.~\ref{eq:model_def1}
  and \ref{eq:model_def2} is well defined. More formally:
  \begin{proposition*}[\ref{prop:well_defined}]

    For every
    realization of the random variables $X_v(0), v\in V$ and for
    all $t\geq1$, $X_v(t)$ is uniquely determined.
  \end{proposition*}

\item
  The process can be calculated efficiently:
\begin{proposition*}[\ref{prop:efficient}]
  It is possible to calculate
  $\{X_v(t) | v\in V\}$, given $\{X_v(t-1) | v \in V\}$, by performing $n$ operations
  of finding the point of an $n$ dimensional
  affine space (as specified by a generating set of size at most $n$) with minimal $L_2$ norm.
  \end{proposition*}

  Calculating the latter is a classical convex optimization problem. See, e.g.,
  \cite{BoydVandenberghe:04}.

\item
  For transitive graphs the process always converges to the optimal estimator:
\begin{proposition*}[\ref{prop:transitive}]
  Let $G$ be a transitive graph (defined below).
  Then $X_v(t)$ converges to $X(\infty)=Z(\infty)$.
\end{proposition*}

\item For graphs of large maximal degree, $X_v(t)$ converge to
  $\mu$:
  \begin{proposition*}[\ref{prop:large_degree}]
    Let $G_n=(V_n,E_n)$ be a family of graphs
    where $|V_n|=n$ and $\max_{v\in V_n}{d_v}\to \infty$. Then
\begin{equation*}
  \lim_{n\to\infty}\sup_{v\in V_n}\lim_{t\to\infty}{\E[\left(X_v(t)-\mu\right)^2]}=0.
\end{equation*}
  \end{proposition*}
  Note by comparison that for any graph,
  \begin{equation*}
  \E[\left(Y(\infty)-\mu\right)^2]=\fourth|E|^{-2}\sum_{v\in V}{d_v^2}.
  \end{equation*}
  In particular for a star on $n$ vertices, as $n\to \infty$ it holds that
  $X(\infty)$ converges to $\mu$ but $Y(\infty)$ does not.

\item
  Finally, for some graphs, the process converges to a limit different than
  $Z(\infty)$ and $Y(\infty)$.
  \begin{theorem*}[\ref{theorem:interval}]
    Let $G=(V,E)$ be the interval of length four where $V = \{a,b,c,d\}$ and
$E = \left\{\{a,b\},\{b,c\},\{c,d\}\right\}$. Then $X_v(t)$ converges to a limit $X(\infty)$,
where
\[
X(\infty) = \fourth\left[(1-\xi)(X_a+X_d) + (1+\xi)(X_b+X_c)\right],
  \quad
\Var[X(\infty)] = \xi,
\]
with $\xi=2-\sqrt{3}=\fourth(1+\sqrt{49}-\sqrt{48})$.
  \end{theorem*}
Note that for this graph
\[
Y(\infty)={1\over 4}\left((1-\fifth)(X_a+X_d)+(1+\fifth)(X_b+X_c)\right),
  \quad\Var[Y(\infty)]=0.26
,
\]
and
\[
 Z(\infty)={1\over 4}\left(X_a+X_b+X_c+X_d\right),
  \quad
\Var[Z(\infty)]=\fourth.
\]

\end{itemize}

\subsubsection{Conjectures}
Showing some supporting results, we conjecture that the process always
converges, and in particular to a state where all agents have the
same estimator.

We present a number of additional open problems and conjectures in the
conclusion.



\section{General Proofs}
\subsection{Process is Well Defined}
\begin{proposition}
  \label{prop:well_defined}
  For every realization of the random variables $X_v(0), v\in V$ and for
  all $t\geq 1$, $X_v(t)$ is uniquely determined.
\end{proposition}
\begin{proof}
  Let $X_v^{(1)}(t)=A$ and $X_v^{(2)}(t)=B$ be minimum variance
  estimators satisfying Eqs.~\ref{eq:model_def1} and \ref{eq:model_def2},
  with variance $V$. Then their average must have variance at least $V$,
  since it also is a linear combination of the estimators from which
  $A$ and $B$ were constructed:
  \begin{eqnarray*}
    V&\leq& \Var[\half(A+B)]\\
    V&\leq& \fourth \Var A + \fourth \Var B +
             \half\Cov(A,B)\\
    V&\leq& \half V + \half\Cov(A,B)\\
    V&\leq& \Cov(A,B).
  \end{eqnarray*}

Since $\Cov(A,B)\leq \sqrt{\Var A\Var B}=V$,
then $\Cov(A,B)=V$ and  $A=B$. Therefore, there exists a unique minimum variance
unbiased estimator and the process is well defined.
\end{proof}

\subsection{The Algorithm for Calculating the Estimator}

We present an efficient algorithm to calculate $X_v(t)$.
Let $E_v(t)=\{X_w(t)|w\in N(v)\}$ be the estimators of agent $v$'s
neighbors at time
$t$. Let  $\bf{C}$ be the covariance matrix of $E_v(t)$, so that
$C_{wu}=\Cov(X_w(t),X_u(t))$. For each $w$, let $x_w$ be a realization
of $X_w(t)$.
Then the log likelihood of $y\in \R$ is
\begin{equation*}
  \log \mathcal{L}(y)=-\sum_{wu}(x_w-y)C^{-1}_{wu}(x_u-y) +\mbox{const},
\end{equation*}
where ${\bf C}^{-1}$ is ${\bf C}$'s pseudo-inverse.
this expression is maximal for
\begin{equation*}
  y = {\sum_{wu}C^{-1}_{wu}x_w\over \sum_{wu}C^{-1}_{wu}}.
\end{equation*}

Hence, the MLE, and therefore also $X_w(t+1)$, equals
\begin{equation}
  \label{eq:algorithm}
  X_{ML}=X_v(t+1)={\sum_{w,u\in N(v)}C^{-1}_{wu}X_w(t)\over \sum_{w,u\in N(v)}C^{-1}_{wu}}.
\end{equation}
Note that $X_v(t+1)$ is also, among all the unbiased estimators of $\mu$
constructible over the estimators in $E_V(t)$, the one with
the minimum variance.

Given this last observation, there exists a simple geometric interpretation
for Eq.~\ref{eq:algorithm}:
\begin{proposition}
  \label{prop:efficient}
  It is possible to calculate
  $\{X_v(t) | v\in V\}$, given $\{X_v(t-1) | v \in V\}$, by performing $n$
  operations of finding the point of an $n$ dimensional
  affine space (as specified by a generating set of size at most $n$)
  with minimal $L_2$ norm.
\end{proposition}
\begin{proof}
Consider an $n$-dimensional vector space $\mathcal{V}$ over $\R$, with an inner product
$\left<\cdot,\cdot\right>$.
Let ${\bf z}$ be some non-zero vector in $\mathcal{V}$, and let $\mathcal{A}\subset\mathcal{V}$ be the affine space
defined by
$\mathcal{A}=\{{\bf x}\in \mathcal{V}|\left<{\bf x},{\bf z}\right>=1\}$.

(This is a generalization of $\mathcal{V}=\textrm{span}(\{X_v|v\in V\})$,
 $\left<X,Y\right>=\Cov(X,Y)$, ${\bf z}=\sum_{v\in V}X_v$ and $\mathcal{A}$
being the set of unbiased estimators).

Given a set of vectors $E=\left\{{\bf x}_k|k=1,\ldots,K\leq n\right\}$, where ${\bf x}_k\in\mathcal{A}$,
let $C_{kl}=\left<{\bf x}_k,{\bf x}_l\right>$.
Then $\mathcal{A}\cap\textrm{span}(E)$ is also an affine space, and
the  minimum $L_2$ norm vector in
  $\mathcal{A}\cap\textrm{span}(E)$ is
  \begin{equation*}
    {\bf x}_{ML}={\sum_{kl}C^{-1}_{kl}{\bf x}_k\over \sum_{kl}C^{-1}_{kl}},
  \end{equation*}
  where ${\bf C}^{-1}$ is the matrix pseudo-inverse of ${\bf C}$.
This equation is identical to Eq.~\ref{eq:algorithm}.
\end{proof}
Note that if ${\bf C}$ is invertible then its pseudo-inverse is equal to its
inverse. Otherwise, there are many linear combinations of the
vectors in $E$ which are equal to the unique ${\bf x}_{ML}$.
The r\^ole of the pseudo-inverse is to facilitate computation: it provides
the linear combination with least sum of squares of the
coefficients \cite{NumericalRecipes:92}.

\subsection{Convergence}
Denote $V_v(t):=\Var[X_v(t)]$ and $C_{vw}(t):=\Cov(X_v(t),X_w(t))$.
\begin{lemma}
\label{lemma:limiting_variance}
All estimators have the same limiting variance: $\exists\rho_\infty\forall v:\;V_v(t)\to \rho_\infty$
\end{lemma}

\begin{proof}
Since, in every iteration, each agent calculates the minimum variance
unbiased estimator over those of his neighbors and its own, then
the variance of the estimator it calculates must be lower than
that of any other unbiased linear combination:
\begin{equation}
  \label{eq:var_mono_general}
  \sum_{w\in N(v)}{\alpha_w}=1\;\Rightarrow \;
  V_v(t+1)\leq \Var\left[\sum_{w\in N(v)}{\alpha_wX_w(t)}\right].
\end{equation}
In particular, for each neighbor $w$ of $v$
\begin{equation}
  \label{eq:var_mono_neigh}
  V_v(t+1)\leq V_w(t),
\end{equation}
and since each vertex is its own neighbor, then
\begin{equation}
  \label{eq:var_mono}
  V_v(t+1)\leq V_v(t).
\end{equation}
Therefore, since the variance of each agent's estimator is
monotonously decreasing (and positive), it
must converge to some $\rho_v$. Now assume $(v,w)\in E$ and $\rho_v<\rho_w$,
then, at some iteration $t$,
 $V_v(t)<\rho_w\leq V_w(t)$. But then, by Eq.~\ref{eq:var_mono_neigh}, we have
$V_w(t+1)\leq V_v(t)<\rho_w$ - a contradiction.
Therefore, $\rho_v$ must equal $\rho_w$, and since the graph is connected,
all agents must converge to the same variance, $\rho_\infty$.
\end{proof}

\begin{lemma} $\forall v,w:\;C_{vw}(t)\to \rho_\infty$
\label{lemma:getting_closer}
\end{lemma}
\begin{proof}
The previous lemma is a special case of this one, for when $v=w$. Otherwise,
for two neighboring agents $v$ and $w$,
for any $\epsilon$, there exists an iteration $t$ where both $\Var[X_v(t)]<\rho_\infty+\epsilon$ and
$\Var[X_w(t)]<\rho_\infty+\epsilon$. Then:
\begin{equation*}
\Var[\half(X_v(t)+X_w(t))]=
  \fourth\left[V_v(t)+V_w(t)+2C_{vw}(t)\right]
  <\half\left[\rho_\infty+\epsilon+C_{vw}(t)\right]
\end{equation*}
and since Eq.~\ref{eq:var_mono_general} implies
$V_v(t+1)\leq \Var[\half(X_v(t)+X_w(t))]$, then
$\rho_\infty< \half\left[\rho_\infty+\epsilon+C_{vw}(t)\right]$ and $C_{vw}(t)\geq\rho_\infty-\epsilon$. Since $C_{vw}$ is also bounded
from above:
$C_{vw}(t)\leq\sqrt{\Var[X_v(t)]\Var[X_w(t)]}<\rho_\infty+\epsilon$, we have demonstrated
that $C_{vw}(t)\to\rho_\infty$ when $v$ and $w$ are neighbors. This implies that the
correlation between neighbors converges to 1, and therefore, since the
graph is finite, all correlations converge to 1 and all covariances
converge to $\rho_\infty$.
\end{proof}

This last lemma implies that if one agent's estimator converges, then
all others' also converge, to the same limit. Even without convergence,
however, it implies that all the estimators converge to their average:
\begin{equation}
\label{eq:getting_closer}
\lim_{t \to \infty}{\Var\left[X_v(t)-{1\over|V|}\sum_{w\in V}{X_w(t)}\right]}=0.
\end{equation}

The following lemma will be used to conjecture that all the estimators do
converge. It states that an estimator is uncorrelated to the difference
between it and any of the estimators which were used to calculate it.
\begin{lemma}
\label{lemma:ortho_iter}
$\forall w\in N(v):\:\Cov(X_v(t+1),X_v(t+1)-X_w(t))=0$
\end{lemma}
\begin{proof}
We examine the estimators
 $\hat{X}(\beta) = X_v(t+1)(1-\beta)+X_w(t)\beta$, which are also
unbiased estimators of $\mu$, and are linear combinations of the estimators from which
$X_v(t+1)$ was constructed. They should all therefore have higher
variance than $X_v(t+1)$. Since $\hat{X}(\beta=0)=X_v(t+1)$, then
\begin{equation*}
 0={\partial\Var[\hat{X}]\over \partial\beta}\Big |_{\beta=0}.
\end{equation*}
Now:
\begin{eqnarray*}
 0&=&{\partial\Var[\hat{X}]\over \partial\beta}\Big |_{\beta=0}\\
 &=&{\partial\left[(1-\beta)^2\Var[X_v(t+1)]+\beta^2\Var[X_w(t)]+2\beta(1-\beta)\Cov(X_v(t+1),X_w(t))\right]\over \partial\beta}\Big |_{\beta=0}\\
 &=&\left[-2(1-\beta)\Var[X_v(t+1)]+2\beta\Var[X_w(t)]+2(1-2\beta)\Cov(X_v(t+1),X_w(t))\right]\Big |_{\beta=0}\\
 &=&-2\Var[X_v(t+1)]+2\Cov(X_v(t+1),X_w(t))\\
 &=&-2\Cov\left(X_v(t+1),X_v(t+1)-X_w(t)\right)\\
\end{eqnarray*}
and so $\Cov\left(X_v(t+1),X_v(t+1)-X_w(t)\right)=0$.
\end{proof}

Note that this implies that $\Var[X_v(t+1)]=\Cov(X_v(t+1),X_w(t))$.

\begin{conjecture}{\bf (Convergence)} $\exists X(\infty)\forall v:\;X_v(t)\to X(\infty)$.
  \label{conj:convergence}
\end{conjecture}

The following observation supports this conjecture:

\begin{eqnarray*}
\Var\left[X_v(t+1)-X_v(t)\right]
&=&\Cov\left(X_v(t+1)-X_v(t),X_v(t+1)-X_v(t)\right)\\
&=&\Cov\left(X_v(t+1),X_v(t+1)-X_v(t)\right)-\Cov\left(X_v(t),X_v(t+1)-X_v(t)\right)
\end{eqnarray*}
Using Lemma~\ref{lemma:ortho_iter}
\begin{eqnarray*}
&=&-\Cov\left(X_v(t),X_v(t+1)-X_v(t)\right)\\
&=&V_v(t)-\Cov(X_v(t+1),X_v(t)),
\end{eqnarray*}
and using it again:
\begin{equation*}
=V_v(t)-V_v(t+1).
\end{equation*}
This implies that if $t_0$ is such that $V_v(t_0) = \rho_\infty+\epsilon$ and therefore
 $\sum_{t= t_0}^\infty{V_v(t)-V_v(t+1)}=\epsilon$, then
\begin{equation*}
  \sum_{t= t_0}^\infty{\Var\left[X_v(t+1)-X_v(t)\right]}=\epsilon.
\end{equation*}

\subsection{Efficiency for Transitive Graphs}
Vertex transitive graphs (henceforth referred to as transitive graphs),
are graphs where all vertices are essentially equivalent,
or ``equally important''. Alternatively, one may say that the graph
``looks the same'' from all vertices. Formally, $G=(V,E)$ is transitive iff,
for every pair of vertices $v, w \in V$ there exists a function $f:V\to V$ which
is a graph automorphism (i.e. $f$ is a bijection and
 $(a,b)\in E\Leftrightarrow (f(a),f(b))\in E$) and maps $v$ to $w$.

\begin{proposition}
  \label{prop:transitive}
  When $G$ is transitive then the process converges and $X(\infty)=Z(\infty)$.
\end{proposition}
\begin{proof}
  By the symmetry of the graph, the average of the agents' estimators
cannot give more weight to one agent's original estimator than to another:
  \begin{equation*}
{1\over |V|}\sum_v{X_v(t)}={1\over |V|}\sum_v{X_v}=Z(\infty),
  \end{equation*}
and hence the average of the agents' estimators is constant and
in particular converges.
By lemma~\ref{lemma:getting_closer} (Eq.~\ref{eq:getting_closer}),
if the average converges then each of the estimators converges to the same
limit:
\begin{equation*}
  \forall v\lim_{t \to \infty}{X_v(t)}=\lim_{t \to \infty}{{1\over |V|}\sum_v{X_v(t)}}=Z(\infty)=X(\infty).
\end{equation*}
\end{proof}

Note that for regular graphs (i.e. graphs where all vertices have the same
degree), which are a superset of transitive graphs, $Y(\infty)=Z(\infty)$.

\section{Analytic Examples}
Complete analytical analysis of these iterations for general graphs
seems difficult, since Eq.~\ref{eq:algorithm} is quadratic. In fact,
we found only two simple examples amenable to complete analysis:
The star, a graph with a central node connected to all other nodes, and
the interval of length four, a graph of four linearly ordered nodes.

In the former, we show that the minimum variance scheme is
efficient, so that $X(\infty)=Z(\infty)$.
In the latter, we show that it isn't,
but that the ``price of anarchy'' is low.

\subsection{High degree graphs and the star}
We consider a graph of $n$ vertices, of which $u$ is the central node
and is connected to all others, and no additional edges exist.

The averaging estimator $Y(\infty)$ gives weight ${n\over 3n-2}$ to $X_u$ and
${2\over 3n-2}$ to the rest. Its variance
is ${n^2+4n-4\over (3n-2)^2}$, which is asymptotically ${1\over 9}$.

On the other hand, $X_u(1)=Z(\infty)$, since node $u$, neighboring all nodes
of the graph, immediately finds the global
minimum variance estimator. In the next iteration, all nodes
$w$ set $X_w(2)=X_u(1)$, and the process essentially halts, since all nodes
have the same estimator, $X_w(2)=X(\infty)=Z(\infty)$, with $\Var[X(\infty)]={1\over n}$.

In general, in graphs of large maximal degree, $X_v(t)$ converges to $\mu$:
\begin{proposition}
  \label{prop:large_degree}
  Let $G_n=(V_n,E_n)$ be a family of graphs
  where $|V_n|=n$ and $\max_{v\in V_n}{d_v}\to \infty$. Then
\begin{equation*}
  \lim_{n\to\infty}\sup_{v\in V_n}\lim_{t\to\infty}{\E[\left(X_v(t)-\mu\right)^2]}=0.
\end{equation*}
\end{proposition}
\begin{proof}
  Since all estimators at all iterations have mean $\mu$, then
  $\E[\left(X_v(t)-\mu\right)^2]=\Var[X_v(t)]$.
  By lemma~\ref{lemma:limiting_variance}, the limiting variances of
  all the agents in a graph $G_n$ are equal to some $\rho_n$, and therefore
  \begin{equation*}
      \lim_{n\to\infty}\sup_{v\in V_n}\lim_{t\to\infty}{\E[\left(X_v(t)-\mu\right)^2]}=0\quad\leftrightarrow\quad\lim_{n\to\infty}\rho_n=0
  \end{equation*}

  The condition $\max_{v\in V_n}{d_v}\to \infty$ implies that
  given $\epsilon>0$, there exists a high enough $N$, so that in any $G_n$ with
  $n>N$ there exists a node $w_n$ with degree $d_{w_n}$ larger than $1/\epsilon$. Then
  $\Var[X_{w_n}(1)]<\epsilon$, since agent $w_n$ would, on the first
  iteration, average the estimators of all its neighbors, resulting
  in a new estimator of variance $1/d_{w_n}$. Since variance never increases
  in the iterative process (lemma~\ref{lemma:limiting_variance}),
  then $\rho_n< \epsilon$ for $n$ larger then some $N$.
  Since this is true for
  arbitrary $\epsilon$, $\rho_n$ goes to zero as $n$ goes to infinity.
\end{proof}

\subsection{Interval of Length four}
We analyze the case of  $G=(V,E)$ where $V = \{a,b,c,d\}$ and
$E = \left\{\{a,b\},\{b,c\},\{c,d\}\right\}$.

\begin{figure}
\begin{center}
\epsfig{file=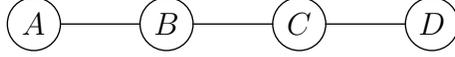}
\caption{Interval of Length Four}
\end{center}
\end{figure}

In Appendix~\ref{app:interval_proof}, we prove
that the process converges with
\begin{equation*}
\Var[X(\infty)]=2-\sqrt{3} =\fourth(1+\sqrt{49}-\sqrt{48}) > \fourth =\Var[Z(\infty)],
\end{equation*}
thus proving that a case exists where $X(\infty)\neq Z(\infty)$.
We also, for this case, derive an asymptotic
convergence rate of $2-\sqrt{3}$.

The averaging estimator $Y(\infty)$ is:
\begin{equation*}
  Y(\infty)=0.2X_a+0.3X_b+0.3X_c+0.2X_d,
\end{equation*}
with
\begin{equation*}
  \Var[Y(\infty)]=0.26.
\end{equation*}
This is slightly lower than $\Var[X(\infty)]$, which equals about $0.268$. However,
the convergence rate (second eigenvalue)
for the averaging process is $\fourth+\sqrt{33}/12\approx 0.73$,
which is significantly slower than the minimum variance process's
asymptotic rate of $2-\sqrt{3}\approx 0.268$.

\section{Numerical Examples and Conjectures}
Numerical simulations on intervals of lengths larger than four
suggest a surprising result.

\subsection{Interval of Arbitrary Length}
Numerical simulations suggest that $X(\infty)$, for intervals of length $n$, approaches
a normal distribution around the center of the interval,
with variance proportional to $n$:
\begin{conjecture}
For interval graphs of length $2n$, index the agents by
 $k\in \{0,\ldots,2n-1\}$.
Then
\begin{equation*}
 X(\infty)=\sum_kA_kX_k
\end{equation*}
where $A_k$ approaches a normal distribution in the sense that
\begin{equation*}
  \lim_{n\to\infty}\sum_k{\left(A_k-C_ne^{-(k-n+1/2)^2/\nu(n)}\right)^2}=0
\quad\mbox{with}\quad\nu(n)\in\Theta(n), \:C_n\in\R.
\end{equation*}
\end{conjecture}
\begin{figure}[htp]
\begin{center}
\epsfig{file=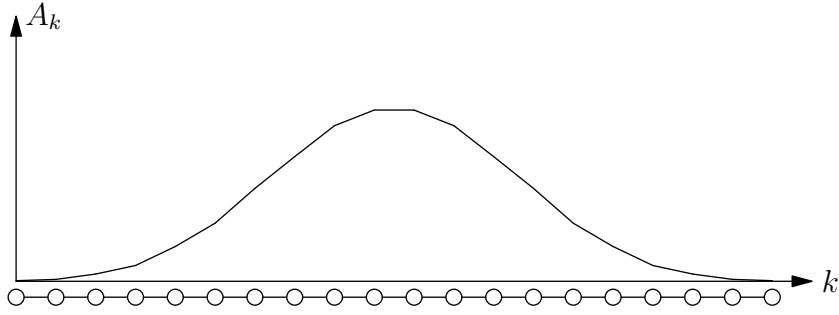}
\caption{Simulation of interval of length 20.}
\end{center}
\end{figure}

 This implies that while $\lim_{n\to\infty}\E[(X(\infty)-\mu)^2]=0$,
$X(\infty)$ quickly becomes less efficient when compared to $Z(\infty)$:
\begin{conjecture}
$\Var[X(\infty)] \propto \sqrt{n}\Var[Z(\infty)]$.
\end{conjecture}

Note that $\Var[Y(\infty)]$, on the other hand,
approaches $\Var[Z(\infty)]$ as $n$ increases, for
intervals of length $n$.

\subsection{Agents with Memory}

A model which is perhaps more natural than the memoryless model is the model
in which the agents remember all their own values from the previous iterations.

\begin{proposition}
If the agents have memory, then the process converges to $Z(\infty)$.
\end{proposition}
\begin{proof}
  Since each vertex $v$ always remembers $X_v=X_v(0)$, then $X_v$ is always
part of the set over which $X_v(t)$ was constructed. Then,
by Lemma~\ref{lemma:ortho_iter}:
\begin{equation}
  \Cov(X_v(t),X_v)=\Var[X_v(t)],
\end{equation}
and by Lemma~\ref{lemma:getting_closer}:
\begin{equation*}
 \forall v,w\in V:\: \lim_{t\to\infty}{\Cov(X_w(t),X_v)}=\lim_{t\to\infty}{\Cov(X_v(t),X_v)}=\lim_{t\to\infty}{\Var[X_v(t)]}=\rho_\infty.
\end{equation*}

This means that for any agent $w$, the covariance of its limit estimator
with each of the original estimators $X_v$ is identical, and so it must
be their average: $X(\infty)=Z(\infty)$.
\end{proof}

This proof relied only on the agents' memory of their original estimators.
Since they also gain more estimators over the iterations, and seemingly
expand the space that they span, we conjecture that:
\begin{conjecture}
  $\forall v\in V:\:X_v(t)=X(\infty)$, for $t\geq |V|$.
\end{conjecture}

\section{Conclusion}
An number of interesting open problems can be raised with respect to this
 model, some of which we conjecture about above:
\begin{itemize}
\item Does it always converge? We conjecture above that this is indeed the case.
\item For what graphs does it converge to the optimal estimator $Z(\infty)$?
\item Otherwise, what is the ``price of anarchy'',
 $\Var[X(\infty)]/\Var[Z(\infty)]$? Is it bounded? We conjecture
  above that it isn't.
\item What is the convergence rate?
\end{itemize}

\appendix
\section{Analysis of Interval of Length Four}

  \begin{theorem}
    \label{theorem:interval}
    Let $G=(V,E)$ be the interval of length four where $V = \{a,b,c,d\}$ and
$E = \left\{\{a,b\},\{b,c\},\{c,d\}\right\}$. Then $X_v(t)$ converges to a limit $X(\infty)$,
where
\[
X(\infty) = \fourth\left[(1-\xi)(X_a+X_d) + (1+\xi)(X_b+X_c)\right],
  \quad
\Var[X(\infty)] = \xi,
\]
with $\xi=2-\sqrt{3}=\fourth(1+\sqrt{49}-\sqrt{48})$.
  \end{theorem}
\begin{proof}
\label{app:interval_proof}
We define $M_{vw}(t) = \Cov(X_v,X_w(t))$, so that each column of ${\bf M}$ is
the coordinates of an agent's estimator at time $t$, viewed as
a vector in the space spanned by $\{X_a, X_b, X_c, X_d\}$. We define
$Z(\infty)$-subtracted $M$ as
$\tilde{M}_{vw}(t)=\Cov(X_v-Z(\infty),X_w(t)-Z(\infty))$, where $Z(\infty)=\fourth(X_a+X_b+X_c+X_d)$, and likewise
define the $Z(\infty)$-subtracted covariance matrix
$\tilde{C}_{vw}(t)=\Cov(X_v(t)-Z(\infty),X_w(t)-Z(\infty))$.

We now shift to an alternative orthonormal basis $B$:

\begin{equation*}
B=
\left\{
\begin{pmatrix}
1/2 \\ 1/2 \\ 1/2 \\ 1/2
\end{pmatrix}(=b_1),
\begin{pmatrix}
-1/\sqrt{2} \\ 0 \\ 0 \\ 1/\sqrt{2}
\end{pmatrix}(=b_2),
\begin{pmatrix}
0 \\ -1/\sqrt{2} \\ 1/\sqrt{2} \\ 0
\end{pmatrix}(=b_3),
\begin{pmatrix}
-1/2 \\ 1/2 \\ 1/2 \\ -1/2
\end{pmatrix}(=b_4)
\right\}
\end{equation*}
The vector $b_1(=2Z(\infty))$ was chosen because its coordinate is
one half in every unbiased estimator (and zero for any $Z(\infty)$-subtracted
unbiased estimator). $b_2$ and $b_3$ are anti-symmetric to
inversion of the interval, a transformation which should leave $X(\infty)$
invariant by the symmetry of the graph.
Therefore we expect their coordinates in $X(\infty)$ to vanish.
We have no freedom, then,
in choosing the last vector, and expect $X(\infty)$ to equal
$\half b_1$ plus some constant $\xi$ times $\half b_4$:
\begin{equation*}
X(\infty) = \fourth\left[(1-\xi)(X_a+X_d) + (1+\xi)(X_b+X_c)\right].
\end{equation*}

Performing two iterations of the process reveals that under this basis,
the  $Z(\infty)$-subtracted coordinates matrix of the
estimators at iteration two, $\tilde{\bf M}_V(2)$, is of the form:
\begin{equation*}
\tilde{\bf M}_B(2) =
\begin{pmatrix}
{0}     & {0}           & {0}           & {0}\\
{x}     & {0}           & {0}           & {-x}\\
{0}     & {z}           & {-z}          & {0}\\
{y}     & {w}           & {w}           & {y}
\end{pmatrix},
\end{equation*}
with
\begin{equation}
  \label{eq:a_copies}
 yw=z^2+w^2.
\end{equation}
Application of another iteration yields a matrix of the same form:
\begin{equation}
\label{eq:matrix_iteration}
\tilde{\bf M}_B(3) =
\begin{pmatrix}
{0}     & {0}           & {0}           & {0}\\
{0}     & {{-xz \over x^2+(y-w)^2}z} & {{xz \over x^2+(y-w)^2}z}      & {0}\\
{z}     & {0}           & {0}           & {-z}\\
{w}     & {{x^2 \over x^2+(y-w)^2}w}       & {{x^2 \over x^2+(y-w)^2}w}   & {w}
\end{pmatrix}
,
\end{equation}
with the relation of Eq.~\ref{eq:a_copies} preserved.

Since the result is a matrix of essentially the same form,
equivalent equations apply for consecutive iterations,
and we may denote as $x_t$, $y_t$, $w_t$ and $z_t$ the corresponding matrix
entries at time $t$.

Since Eq.~\ref{eq:matrix_iteration} implies that
 $y_t = w_{t-1}$ and $x_t = z_{t-1}$, then if $w_t$ and $z_t$ converge then
the process converges and $X(\infty)$ exists. Also:
\begin{equation}  \label{eq:ws}
  w_{t+1} ={z_{t-1}^2 \over z_{t-1}^2+(w_{t-1}-w_t)^2}w_t
\end{equation}
and
\begin{equation}  \label{eq:zeds}
  z_{t+1} ={z_{t-1}z_t \over z_{t-1}^2+(w_{t-1}-w_t)^2}z_t.
\end{equation}
Dividing Eq.~\ref{eq:ws} by Eq.~\ref {eq:zeds}, we discover that:
\begin{equation*}
  {w_{t+1} \over z_{t+1} } = {z_{t-1} \over z_t}{w_t \over z_t},
\end{equation*}
and therefore, by repeated application:
\begin{equation*}
  w_t  = {w_2z_2 \over z_3}{z_t \over z_{t-1}} = {1\over 2}{z_t \over z_{t-1}}.
\end{equation*}
Eq.~\ref{eq:a_copies} can alternatively be written as:
$w_{t-1}w_t=w_t^2+z_t^2$. Then:
\begin{equation*}
w_{t-1} - w_t = z_t^2/w_t=2z_tz_{t-1},
\end{equation*}
and we can write Eq.~\ref{eq:zeds} as:
\begin{equation*}
z_{t+1}={z_{t-1}z_t \over z_{t-1}^2+4z_t^2z_{t-1}^2}z_t
\end{equation*}
or
\begin{equation}
\label{eq:simpzeds}{z_t\over z_{t+1}}={z_{t-1}\over z_t}(1+4z_t^2).
\end{equation}

To solve this recursion we make the following guess:
\begin{equation}
\label{eq:guess}
{z_t\over z_{t+1}}=2+\sqrt{3-4z_t^2},
\end{equation}
which is a solution of the following quadratic equation in $z_t/z_{t+1}$:
\begin{equation*}
{z_t^2\over z_{t+1}^2}-4{z_t\over z_{t+1}}+1+4z_t^2=0.
\end{equation*}
This is equivalent to the following relation:
\begin{equation*}
\Var[X_b(t)]=w_t^2+z_t^2+\fourth=2w_t,
\end{equation*}
upon which we serendipitously stumbled during our examination of this problem.

This guess satisfies Eq.~\ref{eq:simpzeds}, as some manipulation of the
two equations will show. Since $z_2$ and $z_3$ satisfy
Eq.~\ref{eq:guess}, then the rest of the $z$'s must, too.

Since Eq.~\ref{eq:simpzeds}
implies $z_t\to 0$, we can conclude from
${z_{t-1}\over z_t}=2+\sqrt{3-4z_{t-1}^2}$ that
\begin{equation*}
\lim_{t\to\infty}{{z_{t+1}\over z_t}}= 2-\sqrt{3}:=\xi.
\end{equation*}
This is the process's asymptotic convergence rate. Since $w_t={1\over 2}{z_t\over z_{t-1}}$,
then $w_t\to {1\over 2}\xi$, and

\begin{equation*}
X(\infty) = \fourth\left[(1-\xi)(X_a+X_d) + (1+\xi)(X_b+X_c)\right],
\end{equation*}
with
\begin{equation*}
\Var[X(\infty)]= \left({1\over 2}\right)^2 + \left({1\over 2}\xi\right)^2 = \xi.
\end{equation*}
\end{proof}
\bibliographystyle{abbrv} \bibliography{all}

\end{document}